\documentclass[11pt,a4paper]{article}

\usepackage{amsmath,amssymb}
\newtheorem{theorem}{Theorem}[section]
\newtheorem{corollary}{Corollary}[section]

\newtheorem{proposition}{Proposition}[section]
\newtheorem{definition}{Definition}[section]
\newtheorem{remark}{Remark}[section]
\numberwithin{equation}{section}
\newtheorem{example}{Example}[section]
\newenvironment{proof}[1][Proof]{\noindent\textbf{#1.} }{\hfill {$\Box$}}
\numberwithin{equation}{section}

\begin{document}

\title{\textsc{Exponential dichotomies for linear discrete-time systems in Banach
spaces}}
\author{\textsc{Ioan-Lucian Popa \ \ \ Mihail Megan\ \ \ Traian Ceau\c su}}
\date{}
\maketitle

{\footnotesize \noindent \textbf{Abstract.} In this paper we
investigate some dichotomy concepts  for linear difference
equations in Banach spaces. We motivate our approach by
illustrative examples.}

{\footnotesize \vspace{3mm} }


{\footnotesize \vspace{2mm} }

\section{Introduction}
In the mathematical literature of the last decades the asymptotic
behavior of difference equations is one of the most important
subjects due to large applications area (see
\cite{agarwal},\cite{bareira},\cite{dale},\cite{Elaydi1},\cite{mil},\cite{lasmikantham},\cite{massera}).

The classical paper of Perron \cite{perron} served as a starting
point for numerous works on the stability theory. A discrete
variant of Perron's results was given by Ta Li in \cite{tali}.
Several results about exponential dichotomy were obtained by C.V.
Coffman and J.J. Sch$\ddot{a}$ffer \cite{coffman}, J. Kurzweil and
G. Papaschinopoulos \cite{kurzweil}, M. Pinto \cite{pinto}, Y.
Latushkin et al. \cite{latushkin}, P. H. Ngoc and  T. Naito
\cite{ngoc}. Connections between admissibility and uniform
exponential dichotomy of discrete evolution families are given in
\cite{sasu}.

In their notable contribution \cite{bareira} L. Barreira and C.
Valls obtain results in the case
 of a notion of nonuniform exponential dichotomy, which is motivated by ergodic theory. A principal
 motivation for weakening the assumption of uniform exponential behavior is that from the point of view
 of ergodic theory, almost all linear variational equations in a finite dimensional space admit a nonuniform exponential
 dichotomy. Characterizations in terms of Lyapunov function for
 nonuniform exponential dichotomy are given in \cite{bareira1}.

In this paper we study some dichotomy concepts for non-autonomous
linear discrete-time systems in Banach spaces. The most common
classes of exponential dichotomy used in the qualitative theory of
difference equations are the uniform and nonuniform exponential
dichotomy. The present paper considers three concepts of
nonuniform exponential dichotomy and the classical property of
uniform exponential dichotomy for difference equations in Banach
spaces.

The aim of this paper is to define and to exemplify these concepts
and to emphasize connections between them.

The obtained results extend the framework to the study of
dichotomy of difference equations, hold without any requirement on
the coefficients and are applicable to all systems of difference
equations.
\section{Notations. Definitions.}
Let $X$ be a real or complex Banach space and $\mathcal{B}(X)$ the
Banach algebra of all bounded linear operators from $X$ into
itself. The norms of both these spaces will be denoted by
$\parallel .
\parallel.$ Let $\mathbb{N}$ be the set of all positive
integers and $\Delta$ be the set of all pairs $(m,n)$ of positive
integers satisfying the inequality $m \geq n.$ We also denote by
$T$ the set of all triplets $(m,n,p)$ of positive integers with
$(m,n)$ and $(n,p)\in\Delta .$

In this paper we consider linear discrete-time systems of the
 form

\begin{equation*}\tag{$\mathfrak{A}$}\label{A}
x_{n+1}={A}{(n)}x_{n}.
\end{equation*}
where $A:\mathbb{N}\longrightarrow\mathcal{B}(X)$ is a sequence in
$\mathcal{B}(X)$. Then every solution $x=(x_{n})$ of the system
(\ref{A}) is given by
\begin{equation*}
x_{m}=\mathcal{A}(m,n)x_{n},\;\;\;\text{for
all}\;\;\;(m,n)\in\Delta,
\end{equation*}
where $A:\Delta\longrightarrow\mathcal{B}(X)$ is defined by
\begin{equation}\label{eqAmn}
\mathcal{A}(m,n)=\left\{\begin{array}{ll}
{A}{(m)}\cdot \ldots\cdot  {A}{(n+1)},\;\; \;\; m \geq n+1\\
\qquad\quad I\qquad\qquad\quad,\;\;\;\;\; m=n.
\end{array}\right.
\end{equation}
where $I$ is the identity operator on $X$.  For the particular
case when ${A}(n)={A}\in\mathcal{B}(X)$ we  have that
$\mathcal{A}(m,n)={A}^{m-n}$ for all $(m,n)\in\Delta.$  It is
obvious that
\begin{equation*}
\mathcal{A}(m,n)\mathcal{A}(n,p)=\mathcal{A}(m,p),\;\;\;\text{for
all}\;\;\;(m,n,p)\in T.
\end{equation*}

\begin{definition}\label{D:proiector}
An application $P:\mathbb{N}\rightarrow\mathcal{B}(X)$ is said to
be a {\it family of projections} on $X$ if
\begin{equation*}
P^{2}(n)=P(n),
\end{equation*}
for every $n\in\mathbb{N}.$
\end{definition}
\begin{remark} If $P:\mathbb{N}\rightarrow\mathcal{B}(X)$ is a
family of projections then
$Q:\mathbb{N}\rightarrow\mathcal{B}(X),$ with $Q(n)=I-P(n)$ is
also a family of projections on $X,$ and it is called  {\it the
complementary projection} of $P.$ It is obvious that
\begin{equation*}
P(n)Q(n)=Q(n)P(n)=0,
\end{equation*}
for each $n\in\mathbb{N}.$
\end{remark}
\begin{definition}
A family of projections $P:\mathbb{N}\rightarrow\mathcal{B}(X)$ is
said to be {\it compatible with  the system} (\ref{A}), if
\begin{equation*}
{A}(n+1)P(n)=P(n+1){A}(n+1),
\end{equation*}
for every $n\in\mathbb{N}.$
\end{definition}
\begin{remark}
For the particular case when (\ref{A}) is autonomous, i.e.
${A}(n)={A}\in\mathcal{B}(X)$ for all $n\in\mathbb{N},$ and
$P(n)=P$ then $P$ is compatible with system (\ref{A}) if and only
if ${A}P=P{A}.$
\end{remark}
If $P:\mathbb{N}\rightarrow\mathcal{B}(X)$ is a family of
projections compatible with system (\ref{A}) then
\begin{equation*}
\mathcal{A}(m,n)P(n)=P(m)\mathcal{A}(m,n),
\end{equation*}
for all $(m,n)\in\Delta.$

In what follows, we will denote by
$\mathcal{A}_{P}:\Delta\rightarrow\mathcal{B}(X)$ and
$\mathcal{A}_{Q}:\Delta\rightarrow\mathcal{B}(X)$ the mappings
defined by

\begin{equation*}
\mathcal{A}_{P}(m,n)=\left\{\begin{array}{ll}
{A}{(m)}\cdot \ldots\cdot  {A}{(n+1)}P(n),\;\; \;\; m >n\\
\qquad\quad P(n)\qquad\qquad\quad,\;\;\;\;\; m=n.
\end{array}\right.
\end{equation*}
and
\begin{equation*}
\mathcal{A}_{Q}{(m,n)}=\left\{\begin{array}{ll}
{A}{(m)}\cdot \ldots\cdot  {A}{(n+1)}Q(n),\;\; \;\; m >n\\
\qquad\quad Q(n)\qquad\qquad\quad,\;\;\;\;\; m=n.
\end{array}\right.
\end{equation*}
for all $(m,n)\in\Delta.$
\section{Uniform exponential dichotomy}
Let $P:\mathbb{N}\rightarrow\mathcal{B}(X)$ be a family of
projections compatible with system (\ref{A}), thus
%
%
\begin{definition}\label{D:ued}
The linear discrete-time system (\ref{A}) is said to be {\it
P-uniformly exponentially dichotomic} (and denote P-u.e.d.) if
there exist two constants $N \geq 1$ and $\alpha >0$  such that
\begin{eqnarray}\label{eq D:ued}
&e^{\alpha (m-n)}\left(\parallel
{A}{(m)}\cdot\ldots\cdot{A}{(n+1)}P(n)x\parallel
+\parallel Q(n)x\parallel\right)\leq\nonumber\\
 &\leq N
\left(\parallel P(n)x\parallel + \parallel
{A}{(m)}\cdot\ldots\cdot{A}{(n+1)}Q(n)x\parallel\right)
\end{eqnarray}
for all $(m,n,x)\in \Delta\times X.$
\end{definition}
\begin{remark}\label{R:ued}
For every linear discrete-time system (\ref{A}) the following
statements are equivalent:
\begin{enumerate}

\item[{\it{i)}}] (\ref{A}) is P-uniformly exponentially
dichotomic;

\item[{\it{ii)}}] there exist two constants $N \geq 1$ and $\alpha
>0$  such that
\begin{eqnarray*}\label{eq1 R:ued}
&e^{\alpha (m-n)}\left(\parallel {\mathcal{A}_{P}(m,n)}x\parallel
+ \parallel
Q(n)x\parallel\right) \leq\\
&\leq N \left( \parallel P(n)x\parallel +
\parallel {\mathcal{A}_{Q}(m,n)}x\parallel \right)
\end{eqnarray*}
for all $(m,n,x)\in \Delta\times X.$

\item[{\it{iii)}}] there exist two constants $N \geq 1$ and
$\alpha
>0$  such that
\begin{eqnarray*}\label{eq2 R:ued}
&e^{\alpha (m-n)}\left(\parallel {\mathcal{A}_{P}(m,p)}x\parallel
+
\parallel \mathcal{A}_{Q}(n,p)x\parallel\right)\leq\\
&\leq N \left( \parallel \mathcal{A}_{P}(n,p)x\parallel +
\parallel {\mathcal{A}_{Q}(m,p)}x\parallel \right)
\end{eqnarray*}
for all $(m,n,p,x)\in T\times X.$
\end{enumerate}
\end{remark}
\begin{remark}\label{R:aued}
For the particular case when (\ref{A}) is autonomous, i.e.
${A}(n)={A}\in\mathcal{B}(X)$ and $P(n)=P=P^{2}\in\mathcal{B}(X)$
for all $n\in\mathbb{N}$ we have that system (\ref{A}) is
P-uniformly exponentially dichotomic if and only if there exist
two constants $N \geq 1$ and $\alpha >0$ such that
\begin{eqnarray*}\label{eq R:aued}
e^{\alpha (m-n)}\left(\parallel {A}^{m-n}Px\parallel +
\parallel Qx\parallel\right) \leq N \left(
\parallel Px\parallel  + \parallel
{A}^{m-n}Qx\parallel\right)
\end{eqnarray*}
for all $(m,n,x)\in \Delta\times X.$
\end{remark}
The following example presents a system that is P-uniformly
exponentially dichotomic.
\begin{example}\label{E: 1Pued}
Let $X=\mathbb{R}^{2}$ and
${A}:\mathbb{N}\rightarrow\mathcal{B}(\mathbb{R}^{2})$ defined by
\begin{equation*}
{A}(n)(x_{1},x_{2})=\left(\frac{x_{1}}{a_{n}},a_{n}x_{2}\right)
\end{equation*}
for all $(n,x_{1},x_{2})\in\mathbb{N}\times \mathbb{R}^{2},$ where
$a_{n}=e^{n+\frac{1}{2}}.$ Then for
$P:\mathbb{N}\rightarrow\mathcal{B}(\mathbb{R}^{2})$ defined by
\begin{equation*}
P(n)(x_{1},x_{2})=(x_{1},0)
\end{equation*}
for all $(n,x_{1},x_{2})\in\mathbb{N}\times \mathbb{R}^{2},$ we
have that
\begin{equation*}
  \mathcal{A}_{P}(m,n)(x_{1},x_{2}) = \left\{
  \begin{array}{l l}
    \left(e^{\frac{(n+1)^{2}-(m+1)^{2}}{2}}x_{1},0\right) & \quad m>n\\
    \left(x_{1},0\right) & \quad m=n\\
  \end{array} \right. ,
\end{equation*}
\begin{equation*}
  \mathcal{A}_{Q}(m,n)(x_{1},x_{2}) = \left\{
  \begin{array}{l l}
    \left(0,e^{\frac{(m+1)^{2}-(n+1)^{2}}{2}}x_{2}\right) & \quad m>n\\
    \left(0,x_{2}\right) & \quad m=n\\
  \end{array} \right. ,
\end{equation*}
Hence,
\begin{equation*}
e^{\frac{1}{2}(m-n)}\left(\parallel \mathcal{A}_{P}(m,n)x\parallel
+\parallel Q(n)x\parallel\right) \leq \parallel P(n)x\parallel
+\parallel \mathcal{A}_{Q}(m,n)x\parallel,
\end{equation*}
and thus for $N=1$ and $\alpha =\dfrac{1}{2}$ system (\ref{A}) is
P-u.e.d. \hfill $\square$
\end{example}
\section{Nonuniform exponential dichotomy}
%
%
\begin{definition}\label{D:ned}
The linear discrete-time system (\ref{A}) is said to be {\it
P-nonuniformly exponentially dichotomic} (and denote P-n.e.d.) if
there exists a constant $\alpha >0$ and a nondecreasing sequence
of real numbers $N:\mathbb{N}\longrightarrow \mathbb{R}_{+}^{*}$
such that
\begin{eqnarray}\label{eq D:ned}
&e^{\alpha (m-n)}\left(\parallel
{A}{(m)}\cdot\ldots\cdot{A}{(n+1)}P(n)x\parallel
+\parallel Q(n)x\parallel\right)\leq\nonumber\\
 &\leq N(n)\parallel P(n)x\parallel + N(m)\parallel
{A}{(m)}\cdot\ldots\cdot{A}{(n+1)}Q(n)x\parallel
\end{eqnarray}
for all $(m,n,x)\in \Delta\times X.$
\end{definition}
\begin{remark}\label{R:ned}
For every linear discrete-time system (\ref{A}) the following
statements are equivalent:
\begin{enumerate}

\item[{\it{i)}}] (\ref{A}) is P-nonuniformly exponentially
dichotomic;

\item[{\it{ii)}}] there exists a constant $\alpha >0$ and a
nondecreasing sequence of real numbers
$N:\mathbb{N}\longrightarrow \mathbb{R}_{+}^{*}$  such that
\begin{eqnarray}\label{eq1 R:ned}
&e^{\alpha (m-n)}\left(\parallel {\mathcal{A}_{P}(m,n)}x\parallel
+ \parallel
Q(n)x\parallel\right) \leq\nonumber\\
&\leq N(n) \parallel P(n)x\parallel + N(m)\parallel
{\mathcal{A}_{Q}(m,n)}x\parallel
\end{eqnarray}
for all $(m,n,x)\in \Delta\times X.$

\item[{\it{iii)}}] there exists a constant $\alpha >0$ and a
nondecreasing sequence of real numbers
$N:\mathbb{N}\longrightarrow \mathbb{R}_{+}^{*}$  such that
\begin{eqnarray*}\label{eq2 R:ned}
&e^{\alpha (m-n)}\left(\parallel {\mathcal{A}_{P}(m,p)}x\parallel
+
\parallel \mathcal{A}_{Q}(n,p)x\parallel\right)\leq\\
&\leq N(n) \parallel \mathcal{A}_{P}(n,p)x\parallel +N(m)
\parallel {\mathcal{A}_{Q}(m,p)}x\parallel
\end{eqnarray*}
for all $(m,n,p,x)\in T\times X.$
\end{enumerate}
\end{remark}
\begin{remark}\label{R:aned}
For the particular case when (\ref{A}) is autonomous, i.e.
${A}(n)={A}\in\mathcal{B}(X)$ and $P(n)=P=P^{2}\in\mathcal{B}(X)$
for all $n\in\mathbb{N}$ we have that system (\ref{A}) is
P-nonuniformly exponentially dichotomic if and only if there
exists a constant $\alpha >0$ and a nondecreasing sequence  of
real numbers $N:\mathbb{N}\longrightarrow \mathbb{R}_{+}^{*}$ such
that
\begin{eqnarray*}\label{eq R:aned}
e^{\alpha (m-n)}\left(\parallel {A}^{m-n}Px\parallel +
\parallel Qx\parallel\right) \leq N(n)
\parallel Px\parallel  + N(m)\parallel
{A}^{m-n}Qx\parallel
\end{eqnarray*}
for all $(m,n,x)\in \Delta\times X.$
\end{remark}
\begin{remark}
A P-uniformly exponentially dichotomic system (\ref{A}) is
P-nonuniformly exponentially dichotomic. Now, we will present an
example which shows  that the converse implication is not valid.
\end{remark}
\begin{example}\label{E: 1Pned}
Let $X=\mathbb{R}^{2}$ and
${A}:\mathbb{N}\rightarrow\mathcal{B}(\mathbb{R}^{2})$ defined by
\begin{equation*}
{A}(n)(x_{1},x_{2})=\left(b a_{n}{x_{1}}, \dfrac{x_{2}}{b}\right)
\end{equation*}
where
\begin{equation*}
 a_{n} = \left\{
  \begin{array}{l l}
    \dfrac{1}{(n+2)^{c}} & \quad \text{if $n=2k$}\\
    (n+1)^{c} & \quad \text{if $n=2k+1$}\quad \text{for all}\quad (n,x_{1},x_{2})\in\mathbb{N}\times \mathbb{R}^{2},\\
  \end{array} \right.
\end{equation*}
and  $b\in (0,1),$ $c>0.$ Then for
$P:\mathbb{N}\rightarrow\mathcal{B}(\mathbb{R}^{2})$ defined by
\begin{equation*}
P(n)(x_{1},x_{2})=(x_{1},0)
\end{equation*}
for all $(n,x_{1},x_{2})\in\mathbb{N}\times \mathbb{R}^{2},$ we
have that
\begin{equation*}
  \mathcal{A}_{P}(m,n)(x_{1},x_{2}) = \left\{
  \begin{array}{l l}
    \left(b^{m-n}a_{mn}x_{1},0\right) & \quad m>n\\
    (x_{1},0) & \quad m=n\\
  \end{array} \right. ,
\end{equation*}
and
\begin{equation*}
  \mathcal{A}_{Q}(m,n)(x_{1},x_{2}) = \left\{
  \begin{array}{l l}
    \left(0,\dfrac{x_{2}}{b^{m-n}}\right) & \quad m>n\\
    (0,x_{2}) & \quad m=n\\
  \end{array} \right. ,
\end{equation*}
where
\begin{equation*}\label{eq1 E:amn}
  a_{mn} = \left\{
  \begin{array}{l l}
    1 & \quad \text{if $m=2q+1$ and $n=2p+1$}\\
    (n+2)^{c} & \quad \text{if $m=2q+1$ and $n=2p$}\\
    \left(\dfrac{n+2}{m+2}\right)^{c} & \quad \text{if $m=2q$ and $n=2p$}\\
    \left(\dfrac{1}{m+2}\right)^{c} & \quad \text{if $m=2q$ and $n=2p+1$}\\
  \end{array} \right.
\end{equation*}
As we have
\begin{equation*}
b^{m-n}a_{mn}|x_{1}|+|x_{2}|\leq
b^{m-n}N(n)|x_{1}|+b^{m-n}N(m)e^{m-n}|x_{2}|,
\end{equation*}
for all $(m,n,x_{1},x_{2})\in\Delta\times \mathbb{R}^{2}$, where
$N(n)=(n+2)^{c}$ and $\alpha=-\ln b$. By Remark \ref{R:ned} we
obtain that system (\ref{A}) is P-n.e.d.

If we suppose that system (\ref{A}) is P-u.e.d. then there are two
constants $N \geq 1$ and $\alpha >0$ such that
\begin{equation*}
e^{\alpha(m-n)}\left(b^{m-n}a_{mn}|x_{1}|+|x_{2}|\right)\leq
N\left(|x_{1}|+\dfrac{|x_{2}|}{b^{m-n}}\right),
\end{equation*}
for all $(m,n,x_{1},x_{2})\in\Delta\times \mathbb{R}^{2}.$ In
particular, for $m=2q+1$ and $n=2q$ we have that
\begin{equation*}
e^{\alpha}\left(b(2q+2)^{c}|x_{1}|+|x_{2}|\right)\leq
N\left(|x_{1}|+\dfrac{|x_{2}|}{b}\right),
\end{equation*}
for all $q\in\mathbb{N},$ which is a contradiction. Hence, system
(\ref{A}) is not P-u.e.d. \hfill $\square$
\end{example}

\begin{theorem}\label{T: Dned}
The linear discrete-time system (\ref{A}) is P-nonuniformly
exponentially dichotomic if and only if there exists a constant
$d>0$ and a sequence of real numbers $S:\mathbb{N}\longrightarrow
\mathbb{R}_{+}^{*}$ such that
\begin{eqnarray}\label{eq1 T:Dned}
&\sum\limits_{m=n}^{\infty}e^{d(m-n)}\parallel \mathcal{A}_{P}(m,p)x\parallel +\sum\limits_{k=n}^{m}e^{d(m-k)}\parallel \mathcal{A}_{Q}(k,n)x\parallel\leq\nonumber\\
&\leq S(n) \parallel \mathcal{A}_{P}(n,p)x\parallel +S(m)
\parallel {\mathcal{A}_{Q}(m,n)}x\parallel
\end{eqnarray}
for all $(m,n,p,x)\in T\times X.$
\end{theorem}
\begin{proof}
{\it Necessity.} By a simple computation we have that
\begin{eqnarray*}
&\sum\limits_{m=n}^{\infty}e^{d(m-n)}\parallel \mathcal{A}_{P}(m,p)x\parallel +\sum\limits_{k=n}^{m}e^{d(m-k)}\parallel \mathcal{A}_{Q}(k,n)x\parallel\leq\\
&\leq  \frac{e^{\alpha}N(n)}{e^{\alpha}-e^{d}}\parallel
\mathcal{A}_{P}(n,p)x\parallel +
\frac{e^{\alpha}N(m)}{e^{\alpha}-e^{d}}
\parallel {\mathcal{A}_{Q}(m,n)}x\parallel.
\end{eqnarray*}
Hence, for $S(n)=\frac{e^{\alpha}\phi(n)}{e^{\alpha}-e^{d}}$ we
obtain relation (\ref{eq1 T:Dned}), with $\alpha
>0$ and sequence $N(n)$ given by Definition
\ref{D:ned}.

{\it Sufficiency.}  Firstly, we have that
\begin{equation*}
\sum\limits_{m=n}^{\infty}e^{d(m-n)}\parallel
\mathcal{A}_{P}(m,p)x\parallel\leq S(n)\parallel
\mathcal{A}_{P}(n,p)x\parallel
\end{equation*}
and thus
\begin{equation}\label{eq1 DemTDned}
e^{d(m-n)}\parallel \mathcal{A}_{P}(m,n)x\parallel \leq S
(n)\parallel P(n)x\parallel,
\end{equation}
for all $(m,n,x)\in \Delta\times X.$ Similarly, it follows that
\begin{equation}\label{eq2 DemTDned}
e^{d(m-n)}\parallel Q(n)x\parallel \leq S(m)\parallel
\mathcal{A}_{Q}(m,n)x\parallel .
\end{equation}
Finally, using inequalities (\ref{eq1 DemTDned}) and (\ref{eq2
DemTDned}) and Remark \ref{R:ned} we can conclude that system
(\ref{A}) is P-n.e.d.
\end{proof}

As a particular case, we obtain a characterization of P-uniform
exponential dichotomy given by
\begin{corollary}\label{C1: Datkoued}
The linear discrete-time system (\ref{A}) is P- uniformly
exponentially dichotomic if and only if there exist two  constants
$D,d >0$  such that
\begin{eqnarray*}\label{eq1 T:datkoued}
&\sum\limits_{j=m}^{\infty} e^{d(j-n)}\parallel
{\mathcal{A}_{P}(j,n)} x\parallel + \sum\limits_{k=n}^{m}
e^{d(m-k)}\parallel
\mathcal{A}_{Q}(k,n)x\parallel\leq\\
&\leq D \left(\parallel \mathcal{A}_{P}(m,n)x\parallel +
\parallel \mathcal{A}_{Q}{(m,n)}x\parallel\right),
\end{eqnarray*}
for all $(m,n,x)\in \Delta\times X.$
\end{corollary}
Another characterization of the P-uniform exponential dichotomy
property is given by
\begin{corollary}
The linear discrete-time system (\ref{A}) is P- uniformly
exponentially dichotomic if and only if there exists $D
>0$  such that
\begin{eqnarray*}\label{eq21 T:datkoued}
&\sum\limits_{j=m}^{\infty} \parallel {\mathcal{A}_{P}(j,n)}
x\parallel + \sum\limits_{k=n}^{m} \parallel
\mathcal{A}_{Q}(k,n)x\parallel\leq\\
&\leq D \left(\parallel \mathcal{A}_{P}(m,n)x\parallel +
\parallel \mathcal{A}_{Q}{(m,n)}x\parallel\right),
\end{eqnarray*}
for all $(m,n,x)\in \Delta\times X.$
\end{corollary}
\begin{proof}
{\it Necessity.} It is a simple verification as in the proof ot
Theorem \ref{T: Dned}.

{\it Sufficiency.} It is immediate from Corollary \ref{C1:
Datkoued}.
\end{proof}
\begin{remark}
An equivalent variant of the preceding Corollary was proved by P.
Preda and  M. Megan in \cite{megan2}.
\end{remark}
%
\section{Exponential dichotomy}
%
%
\begin{definition}\label{D:ed}
The linear discrete-time system (\ref{A}) is said to be {\it
P-exponentially dichotomic} (and denote P-e.d.) if there exist the
constants $N \geq 1,$  $\alpha >0$ and $\beta \geq 0$  such that
\begin{eqnarray}\label{eq D:ed}
&e^{\alpha (m-n)}\left(\parallel
{A}{(m)}\cdot\ldots\cdot{A}{(n+1)}P(n)x\parallel
+\parallel Q(n)x\parallel\right)\leq\nonumber\\
 &\leq N\left(e^{\beta n}\parallel P(n)x\parallel + e^{\beta m}\parallel
{A}{(m)}\cdot\ldots\cdot{A}{(n+1)}Q(n)x\parallel\right)
\end{eqnarray}
for all $(m,n,x)\in \Delta\times X.$
\end{definition}
\begin{remark}\label{R:ed}
For every linear discrete-time system (\ref{A}) the following
statements are equivalent:
\begin{enumerate}

\item[{\it{i)}}] (\ref{A}) is P-exponentially dichotomic;

\item[{\it{ii)}}] there exist the constants $N \geq 1,$ $\alpha
>0$ and $\beta \geq 0$  such that
\begin{eqnarray*}\label{eq1 R:ed}
&e^{\alpha (m-n)}\left(\parallel {\mathcal{A}_{P}(m,n)}x\parallel
+
\parallel Q(n)x\parallel\right) \leq\\
&\leq  N\left(e^{\beta n}\parallel P(n)x\parallel + e^{\beta m}
\parallel {\mathcal{A}_{Q}(m,n)}x\parallel\right)
\end{eqnarray*}
for all $(m,n,x)\in \Delta\times X.$

\item[{\it{iii)}}] there exist the constants $N \geq 1,$ $\alpha
>0$ and $\beta \geq 0$  such that
\begin{eqnarray*}\label{eq2 R:ed}
&e^{\alpha (m-n)}\left(\parallel {\mathcal{A}_{P}(m,p)}x\parallel
+
\parallel \mathcal{A}_{Q}(n,p)x\parallel\right)\leq\\
&\leq N\left(e^{\beta n} \parallel \mathcal{A}_{P}(n,p)x\parallel
+ e^{\beta m}\parallel {\mathcal{A}_{Q}(m,p)}x\parallel\right)
\end{eqnarray*}
for all $(m,n,p,x)\in T\times X.$
\end{enumerate}
\end{remark}
\begin{remark}\label{R:aed}
For the particular case when (\ref{A}) is autonomous, i.e.
${A}(n)={A}\in\mathcal{B}(X)$ and $P(n)=P=P^{2}\in\mathcal{B}(X)$
for all $n\in\mathbb{N}$ we have that system (\ref{A}) is
P-exponentially dichotomic if and only if there exist the
constants $N \geq 1,$ $\alpha >0$ and $\beta \geq 0$  such that
\begin{eqnarray*}\label{eq R:aed}
e^{\alpha (m-n)}\left(\parallel {A}^{m-n}Px\parallel +
\parallel Qx\parallel\right) \leq N \left(e^{\beta n}\parallel Px\parallel +
e^{\beta m}\parallel {A}^{m-n}Qx\parallel\right)
\end{eqnarray*}
for all $(m,n,x)\in \Delta\times X.$
\end{remark}
\begin{theorem}\label{T: DPed}
The linear discrete-time system (\ref{A}) is P-exponentially
dichotomic if and only if there exist some constants $D \geq 1,$
$d >0$ and $c \geq 0$ such that
\begin{eqnarray}\label{eq T: DPed}
&\sum\limits_{m=n}^{\infty} e^{d(m-n)}\parallel
{\mathcal{A}_{P}(m,p)} x\parallel +\sum\limits_{k=n}^{m}
e^{d(m-k)}\parallel \mathcal{A}_{Q}(k,n)x\parallel
 \leq\nonumber\\
&\leq D\left( e^{c n} \parallel \mathcal{A}_{P}(n,p)x\parallel +
e^{c m}
\parallel \mathcal{A}_{Q}{(m,n)}x\parallel\right)
\end{eqnarray}
for all $(m,n,p,x)\in T\times X.$
\end{theorem}
\begin{proof}
{\it Necessity.} For $d\in (0,\alpha)$ we have that
\begin{eqnarray*}
&\sum\limits_{m=n}^{\infty} e^{d(m-n)}\parallel
{\mathcal{A}_{P}(m,p)} x\parallel +\sum\limits_{k=n}^{m}
e^{d(m-k)}\parallel
\mathcal{A}_{Q}(k,n)x\parallel\leq\\
&\leq N \parallel {\mathcal{A}_{P}(n,p)}x\parallel
\sum\limits_{m=n}^{\infty}
e^{d(m-n)} e^{-\alpha (m-n)}e^{\beta n}+\\
&+N e^{(\beta+d-\alpha)m}\parallel
{\mathcal{A}_{Q}(m,n)}x\parallel \sum\limits_{k=n}^{m}
e^{(\alpha -d)k}=\\
&\leq\frac{N e^{\alpha}}{e^{{\alpha}}-e^{d}} e^{\beta
n}\parallel{\mathcal{A}_{P}(n,p)}x\parallel + \frac{N
e^{\alpha}}{e^{\alpha}-e^{d}} e^{\beta m}
\parallel{\mathcal{A}_{Q}(m,n)}x\parallel.
\end{eqnarray*}
Hence, for $D=1+\frac{N e^{\alpha}}{e^{{\alpha}}-e^{d}}$ and
$c=\beta$  we obtain relation (\ref{eq T: DPed}), with constants
$N,$ $\alpha$ and $\beta$ offered by Definition \ref{D:ed}.

{\it Sufficiency.} According to relation (\ref{eq T: DPed}) we
have that
\begin{eqnarray*}
&\sum\limits_{m=n}^{\infty} e^{d(m-n)}\parallel
{\mathcal{A}_{P}(m,n)}
x\parallel +\sum\limits_{k=n}^{m} e^{d(m-k)}\parallel \mathcal{A}_{Q}(k,n)x\parallel\geq\\
&\geq e^{d(m-n)}\parallel \mathcal{A}_{P}(m,n)x\parallel +
e^{d(m-n)}\parallel {Q}(n)x\parallel.
\end{eqnarray*}
Hence,
\begin{eqnarray*}
&e^{d (m-n)}\left(\parallel {\mathcal{A}_{P}(m,n)}x\parallel +\parallel Q(n)x\parallel\right) \leq\\
&\leq  D\left(e^{c n}\parallel P(n)x\parallel + e^{c m}
\parallel {\mathcal{A}_{Q}(m,n)}x\parallel\right)
\end{eqnarray*}
for all $(m,n,x)\in \Delta\times X.$ Now, using Remark \ref{R:ed}
we obtain that system (\ref{A}) is P-e.d., which completes the
proof.
\end{proof}
\section{Strong exponential dichotomy}
A particular concept of P-exponentially dichotomic system is
defined by
%
%
\begin{definition}\label{D:sed}
The linear discrete-time system (\ref{A}) is said to be {\it
P-strongly exponentially dichotomic} (and denote P-s.e.d.) if
there exist the constants $N \geq 1,$  $\alpha >0$ and $\beta
\in[0,\alpha)$ such that
\begin{eqnarray}\label{eq D:sed}
&e^{\alpha (m-n)}\left(\parallel
{A}{(m)}\cdot\ldots\cdot{A}{(n+1)}P(n)x\parallel
+\parallel Q(n)x\parallel\right)\leq\nonumber\\
 &\leq N\left(e^{\beta n}\parallel P(n)x\parallel + e^{\beta m}\parallel
{A}{(m)}\cdot\ldots\cdot{A}{(n+1)}Q(n)x\parallel\right)
\end{eqnarray}
for all $(m,n,x)\in \Delta\times X.$
\end{definition}
\begin{remark}\label{R:sed}
For every linear discrete-time system (\ref{A}) the following
statements are equivalent:
\begin{enumerate}

\item[{\it{i)}}] (\ref{A}) is P-strongly exponentially dichotomic;

\item[{\it{ii)}}] there exist the constants $N \geq 1,$ $\alpha
>0$ and $\beta \in[0,\alpha)$  such that
\begin{eqnarray*}\label{eq1 R:sed}
&e^{\alpha (m-n)}\left(\parallel {\mathcal{A}_{P}(m,n)}x\parallel
+
\parallel Q(n)x\parallel\right) \leq\\
&\leq  N\left(e^{\beta n}\parallel P(n)x\parallel + e^{\beta m}
\parallel {\mathcal{A}_{Q}(m,n)}x\parallel\right)
\end{eqnarray*}
for all $(m,n,x)\in \Delta\times X.$

\item[{\it{iii)}}] there exist the constants $N \geq 1,$ $\alpha
>0$ and $\beta \in[0,\alpha)$  such that
\begin{eqnarray*}\label{eq2 R:sed}
&e^{\alpha (m-n)}\left(\parallel {\mathcal{A}_{P}(m,p)}x\parallel
+
\parallel \mathcal{A}_{Q}(n,p)x\parallel\right)\leq\\
&\leq N\left(e^{\beta n} \parallel \mathcal{A}_{P}(n,p)x\parallel
+ e^{\beta m}
\parallel {\mathcal{A}_{Q}(m,p)}x\parallel\right)
\end{eqnarray*}
for all $(m,n,p,x)\in T\times X.$
\end{enumerate}
\end{remark}
\begin{remark}\label{R:ased}
For the particular case when (\ref{A}) is autonomous, i.e.
${A}(n)={A}\in\mathcal{B}(X)$ and $P(n)=P=P^{2}\in\mathcal{B}(X)$
for all $n\in\mathbb{N}$ we have that system (\ref{A}) is
P-strongly exponentially dichotomic if and only if there exist the
constants $N \geq 1,$ $\alpha
>0$ and $\beta \in[0,\alpha)$  such that
\begin{eqnarray*}\label{eq R:ased}
e^{\alpha (m-n)}\left(\parallel {A}^{m-n}Px\parallel +
\parallel Qx\parallel\right) \leq N\left(e^{\beta n}\parallel Px\parallel  + e^{\beta m}\parallel
{A}^{m-n}Qx\parallel\right)
\end{eqnarray*}
for all $(m,n,x)\in \Delta\times X.$
\end{remark}
A characterization of the strong exponential dichotomy is given by
\begin{proposition}\label{T: datkosed}
The linear discrete-time system (\ref{A}) is P-strongly
exponentially dichotomic if and only if there exist the constants
$D\geq 1,$ $d >0$ and $0\leq c<d$ such that relation (\ref{eq T:
DPed}) it is true   for all $(m,n,p,x)\in T\times X.$
\end{proposition}
\begin{proof}
It results from Definition \ref{D:sed} and the proof of Theorem
\ref{T: DPed}.
\end{proof}
\begin{remark}
A P-uniformly exponentially dichotomic system (\ref{A}) is
P-strongly exponentially dichotomic. The reciprocal statements is
not true, as shown in what follows.
\end{remark}
\begin{example}\label{E: xsed-ued}
Let $X=\mathbb{R}^{2}$ and
${A}:\mathbb{N}\rightarrow\mathcal{B}(\mathbb{R}^{2})$ defined by
\begin{equation*}
{A}(n)(x_{1},x_{2})=\left(c_{1}a_{n}{x_{1}},c_{2}a_{n}x_{2}\right)
\end{equation*}
for all $(n,x_{1},x_{2})\in\mathbb{N}\times \mathbb{R}^{2},$ where
$c_{1}$ and $c_{2}$ are two positive constants and
\begin{equation*}
 a_{n} = \left\{
  \begin{array}{l l}
    e^{-n} & \quad \text{if $n=2k$}\\
    e^{n+1} & \quad \text{if $n=2k+1$}\\
  \end{array} \right.
\end{equation*}
Then for $P:\mathbb{N}\rightarrow\mathcal{B}(\mathbb{R}^{2})$
defined by
\begin{equation*}
P(n)(x_{1},x_{2})=(x_{1},0)
\end{equation*}
for all $(n,x_{1},x_{2})\in\mathbb{N}\times \mathbb{R}^{2},$ we
have that
\begin{equation*}
  \mathcal{A}_{P}(m,n)(x_{1},x_{2}) = \left\{
  \begin{array}{l l}
    \left(c_{1}^{m-n}a_{mn}x_{1},0\right) & \quad m>n\\
    (x_{1},0) & \quad m=n\\
  \end{array} \right. ,
\end{equation*}
\begin{equation*}
  \mathcal{A}_{Q}(m,n)(x_{1},x_{2}) = \left\{
  \begin{array}{l l}
    \left(0,c_{2}^{m-n}a_{mn}x_{2}\right) & \quad m>n\\
    (0,x_{2}) & \quad m=n\\
  \end{array} \right. ,
\end{equation*}
where
\begin{equation*}\label{eq2 E:amn}
  a_{mn} = \left\{
  \begin{array}{l l}
    1 & \quad \text{if $m=2q$ and $n=2p$}\\
    e^{-n-1} & \quad \text{if $m=2q$ and $n=2p+1$}\\
    e^{m+1} & \quad \text{if $m=2q+1$ and $n=2p$}\\
    e^{m-n} & \quad \text{if $m=2q+1$ and $n=2p+1$}\\
  \end{array} \right. .
\end{equation*}
If we suppose that system (\ref{A}) is P-u.e.d. then there exist
two constants $\alpha >0$ and $N\geq 1$ such that
\begin{equation*}
e^{\alpha (m-n)}\left(c_{1}^{m-n}a_{mn}|x_{1}|+|x_{2}|\right)\leq
N\left(|x_{1}|+c_{2}^{m-n}a_{mn}|x_{2}|\right)
\end{equation*}
for all $(m,n,x_{1},x_{2})\in \Delta\times \mathbb{R}^{2}.$ But
for $m=2k+1,$ $n=2k,$ $x_{1}\neq 0$ and $x_{2}=0$ we have that
\begin{equation*}
e^{\alpha}c_{1}e^{2k+2}\leq N
\end{equation*}
 which is a contradiction.

Further, for $c_{1}=\dfrac{1}{e^{4}},$ $c_{2}=e^{2},$ $\alpha=2,$
$\beta=1$ and $N=e$ we have that
\begin{equation*}
e^{\alpha (m-n)}\left(\parallel \mathcal{A}_{P}(m,n)x\parallel
+\parallel Q(n)x\parallel\right) =e^{\alpha
(m-n)}\left(c_{1}^{m-n}a_{mn}|x_{1}|+|x_{2}|\right)
\end{equation*}
\begin{equation*}
  =\left\{
  \begin{array}{l l}
    e^{\alpha(m-n)}\left(c_{1}^{m-n}|x_{1}|+|x_{2}|\right) & \quad \text{if $m=2q$ and $n=2p$}\\
    e^{\alpha(m-n)}\left(c_{1}^{m-n}e^{-n-1}|x_{1}|+|x_{2}|\right) & \quad \text{if $m=2q$ and $n=2p+1$}\\
    e^{\alpha(m-n)}\left(c_{1}^{m-n}e^{m+1}|x_{1}|+|x_{2}|\right) & \quad \text{if $m=2q+1$ and $n=2p$}\\
    e^{\alpha(m-n)}\left(c_{1}^{m-n}e^{m-n}|x_{1}|+|x_{2}|\right) & \quad \text{if $m=2q+1$ and $n=2p+1$}\\
  \end{array} \right.
\end{equation*}
\begin{equation*}
  \leq\left\{
  \begin{array}{l l}
    e^{-(m-n)}|x_{1}|+c_{2}^{m-n}|x_{2}| & \quad \text{if $m=2q$ and $n=2p$}\\
     e^{-(m+1)}|x_{1}|+c_{2}^{m-n}|x_{2}| & \quad \text{if $m=2q$ and $n=2p+1$}\\
     e^{n+1}|x_{1}|+c_{2}^{m-n}|x_{2}| & \quad \text{if $m=2q+1$ and $n=2p$}\\
     |x_{1}|+c_{2}^{m-n}|x_{2}| & \quad \text{if $m=2q+1$ and $n=2p+1$}\\
  \end{array} \right.
\end{equation*}
 \begin{equation*}
  \leq\left\{
  \begin{array}{l l}
    N\left(e^{\beta n}|x_{1}|+e^{\beta m}c_{2}^{m-n}|x_{2}|\right) & \quad \text{if $m=2q$ and $n=2p$}\\
    N\left(e^{\beta n}|x_{1}|+e^{\beta m}c_{2}^{m-n}e^{-n-1}|x_{2}|\right)  & \quad \text{if $m=2q$ and $n=2p+1$}\\
     N\left(e^{\beta n}|x_{1}|+e^{\beta m}c_{2}^{m-n}e^{m+1}|x_{2}|\right)  & \quad \text{if $m=2q+1$ and $n=2p$}\\
     N\left(e^{\beta n}|x_{1}|+e^{\beta m}c_{2}^{m-n}e^{m-n}|x_{2}|\right)  & \quad \text{if $m=2q+1$ and $n=2p+1$}\\
  \end{array} \right.
\end{equation*}
\begin{equation}\label{eq E:sed}
 =N(e^{\beta n}\parallel P(n)x\parallel +e^{\beta
m}\parallel \mathcal{A}_{Q}(m,n)x\parallel),
\end{equation}
for all $(m,n,x)\in \Delta\times X.$

Hence, the system (\ref{A}) is P-s.e.d., which completes the proof
 \hfill $\square$
\end{example}
\begin{remark}
It is obvious that a P-strongly exponentially dichotomic system
(\ref{A}) is P-exponentially dichotomic. The following example
emphasize the difference between these concepts and presents that
inverse implication is not true.
\end{remark}
\begin{example}\label{E: 4sed}
Let (\ref{A}) be the linear discrete-time system and
$P:\mathbb{N}\rightarrow\mathcal{B}(\mathbb{R}^{2})$ the
projections family considered in Example \ref{E: xsed-ued} with
$c_{1}=e^{-\frac{3}{2}}$ and $c_{2}=e^{\frac{1}{2}}$.

If we suppose that system (\ref{A}) is P-s.e.d. then there exist
the constants $N\geq 1,$ $\alpha
>0,$ $\beta \geq 0$ with $0\leq \beta <\alpha$ such that
\begin{equation*}
e^{\alpha (m-n)}\left(\parallel \mathcal{A}_{P}(m,n)x\parallel
+\parallel Q(n)x\parallel\right) =e^{\alpha
(m-n)}\left(e^{-\frac{3(m-n)}{2}}a_{mn}|x_{1}|+|x_{2}|\right)
\end{equation*}
\begin{equation*}
  =\left\{
  \begin{array}{l l}
    e^{(\alpha-\frac{3}{2})(m-n)}|x_{1}|+e^{\alpha(m-n)}|x_{2}| & \quad \text{if $m=2q$ and $n=2p$}\\
    e^{(\alpha-\frac{3}{2})(m-n)}e^{-n-1}|x_{1}|+e^{\alpha(m-n)}|x_{2}| & \quad \text{if $m=2q$ and $n=2p+1$}\\
    e^{(\alpha-\frac{3}{2})(m-n)}e^{m+1}|x_{1}|+e^{\alpha(m-n)}|x_{2}| & \quad \text{if $m=2q+1$ and $n=2p$}\\
    e^{(\alpha-\frac{3}{2})(m-n)}e^{m-n}|x_{1}|+e^{\alpha(m-n)}|x_{2}| & \quad \text{if $m=2q+1$ and $n=2p+1$}\\
  \end{array} \right.
\end{equation*}
 \begin{equation*}
  \leq\left\{
  \begin{array}{l l}
    N\left(e^{\beta n}|x_{1}|+e^{\beta m}e^{\frac{m-n}{2}}|x_{2}|\right) & \quad \text{if $m=2q$ and $n=2p$}\\
    N\left(e^{\beta n}|x_{1}|+e^{\beta m}e^{\frac{m-n}{2}}e^{-n-1}|x_{2}|\right)  & \quad \text{if $m=2q$ and $n=2p+1$}\\
     N\left(e^{\beta n}|x_{1}|+e^{\beta m}e^{\frac{m-n}{2}}e^{m+1}|x_{2}|\right)  & \quad \text{if $m=2q+1$ and $n=2p$}\\
     N\left(e^{\beta n}|x_{1}|+e^{\beta m}e^{\frac{m-n}{2}}e^{m-n}|x_{2}|\right)  & \quad \text{if $m=2q+1$ and $n=2p+1$}\\
  \end{array} \right.
\end{equation*}
\begin{equation}\label{eq2 E:sed}
 =N(e^{\beta n}\parallel P(n)x\parallel +e^{\beta
m}\parallel \mathcal{A}_{Q}(m,n)x\parallel),
\end{equation}
for all $(m,n,x)\in \Delta\times X.$

In particular, for $x_{1}\neq 0,$ $x_{2}=0,$ $m=2q+1$ and $n=2p$
it follows that
\begin{equation*}
e^{(\alpha -\frac{1}{2})(m-n)}e^{n+1}\leq N e^{\beta n}
\end{equation*}
which is true for $\alpha -\frac{1}{2}\leq 0,$ $N \geq e$ and
$\beta \geq 1.$ It results that $0<\alpha\leq \frac{1}{2} < 1\leq
\beta$ which is a contradiction with the hypotheses  that $0\leq
\beta <\alpha.$ Hence, (\ref{A}) is not P-s.e.d.

Finally, we observe that for $\alpha=\frac{1}{2},$ $\beta=1$ and
$N=e$ relation (\ref{eq2 E:sed}) it is verified and thus the
system (\ref{A}) is P-e.d.
 \hfill $\square$
\end{example}
There exist P-nonuniformly exponentially dichotomic systems that
are not P-exponentially dichotomic, as in the next.
\begin{example}
Let $X=\mathbb{R}^{2}$ and
${A}:\mathbb{N}\rightarrow\mathcal{B}(\mathbb{R}^{2})$ defined by
\begin{equation*}
{A}(n)(x_{1},x_{2})=\left(ca_{n}{x_{1}},\dfrac{x_{2}}{c}\right)
\end{equation*}
for all $(n,x_{1},x_{2})\in\mathbb{N}\times \mathbb{R}^{2},$ where
$c$ is a positive constant  and
\begin{equation*}
 a_{n} = \left\{
  \begin{array}{l l}
    e^{n(1+2^{n})} & \quad \text{if $n=2k$}\\
    e^{-(n+1)(1+2^{n+1})} & \quad \text{if $n=2k+1$}\\
  \end{array} \right.
\end{equation*}
Then for $P:\mathbb{N}\rightarrow\mathcal{B}(\mathbb{R}^{2})$
defined by
\begin{equation*}
P(n)(x_{1},x_{2})=(x_{1},0)
\end{equation*}
for all $(n,x_{1},x_{2})\in\mathbb{N}\times \mathbb{R}^{2},$ we
have that
\begin{equation*}
  \mathcal{A}_{P}(m,n)(x_{1},x_{2}) = \left\{
  \begin{array}{l l}
    \left(c^{m-n}a_{mn}x_{1},0\right) & \quad m>n\\
    (x_{1},0) & \quad m=n\\
  \end{array} \right. ,
\end{equation*}
and
\begin{equation*}
  \mathcal{A}_{Q}(m,n)(x_{1},x_{2}) = \left\{
  \begin{array}{l l}
    \left(0,\dfrac{x_{2}}{c^{m-n}}\right) & \quad m>n\\
    (0,x_{2}) & \quad m=n\\
  \end{array} \right. ,
\end{equation*}
where
\begin{equation*}\label{eq1 E:amn}
  a_{mn} = \left\{
  \begin{array}{l l}
    1 & \quad \text{if $m=2q$ and $n=2p$}\\
    e^{(n+1)(1+2^{n+1})} & \quad \text{if $m=2q$ and $n=2p+1$}\\
     e^{-(m+1)(1+2^{m+1})} & \quad \text{if $m=2q+1$ and $n=2p$}\\
    e^{(n+1)(1+2^{n+1})}e^{-(m+1)(1+2^{m+1})} & \quad \text{if $m=2q+1$ and $n=2p+1$}\\
  \end{array} \right.
\end{equation*}

Let us suppose that system (\ref{A}) is P-e.d. There exist some
constants $N\geq 1,$ $\alpha >0$ and $\beta \geq 0$ such that
\begin{equation*}
e^{\alpha (m-n)}\left(c^{m-n}a_{mn}|x_{1}|+|x_{2}|\right)\leq
N\left( e^{\beta n}|x_{1}| +\dfrac{e^{\beta
m}}{c^{m-n}}|x_{2}|\right)
\end{equation*}
for all $(m,n,x_{1},x_{2})\in \Delta\times \mathbb{R}^{2}.$
Further, if we consider $x_{1}\neq 0$ and $x_{2} =0$ it follows
that
\begin{equation}\label{eqx E:Pned}
\left(e^{\alpha}c\right)^{m-n}a_{mn}\leq N e^{\beta n}.
\end{equation}
There are three different cases at this point.

{\it Case 1.} If $e^{\alpha}c=1$ then $a_{mn}\leq N e^{\beta n}$
which for $n=2p+1$ and $m=n+1$ leads to
\begin{equation*}
\dfrac{n+1}{n}\left(1+2^{n+1}\right)\leq \dfrac{\ln N}{n} +\beta
\end{equation*}
which is false.

{\it Case 2.} If $e^{\alpha}c >1$ then for $m=2q$ and $n=2$
relation (\ref{eqx E:Pned}) became
\begin{equation*}
\left(e^{\alpha c}\right)^{2q-2}\leq N e^{2\beta}
\end{equation*}
which is false.

{\it Case 3.} If $0 < e^{\alpha}c < 1$ then for $n=2p+1$ and
$m=n+1$ relation (\ref{eqx E:Pned}) became
\begin{equation*}
\dfrac{n+1}{n}\left(1+2^{n+1}\right)\leq \dfrac{\ln
(N/e^{\alpha}c)}{n} +\beta
\end{equation*}
which is also false. Hence, system (\ref{A}) is not P-e.d.

Finally, we remark that for $c=\dfrac{1}{e},$
 $\alpha=1$ and $N(n)=e^{(n+1)(1+2^{n+1})}$ the inequality (\ref{eq1 R:ned}) is satisfied and thus by Remark \ref{R:ned} the system (\ref{A})
is P-n.e.d. \hfill $\square$
\end{example}

\section*{Acknowledgement}
$^{\rm b}$This work is financially supported from the Exploratory
Research Grant PN II ID 1080 No. 508/2009 of the Romanian Ministry
of Education, Research and Innovation.
%

\vspace{5mm}

\noindent\begin{tabular}[t]{ll}

Ioan-Lucian Popa \\
Faculty of Mathematics and Computer Science \\
West University of Timi\c{s}oara \\
email: \texttt{popa@math.uvt.ro}

\vspace{3mm}\\

Mihail Megan \\
Faculty of Mathematics and Computer Science \\
West University of Timi\c{s}oara \\
email: \texttt{megan@math.uvt.ro}

\vspace{3mm}\\

Traian Ceau\c su \\
Faculty of Mathematics and Computer Science \\
West University of Timi\c{s}oara \\
email: \texttt{ceausu@math.uvt.ro}
\end{tabular}

\end{document}